\documentclass{amsproc}
\usepackage{graphicx}
\usepackage{amscd}
\usepackage{amsmath}
\usepackage{amsfonts}
\usepackage{amssymb}
\theoremstyle{plain}

\newtheorem{corollary}{Corollary}

\newtheorem{definition}{Definition}
\newtheorem{example}{Example}

\newtheorem{remark}{Remark}

\newtheorem{theorem}{Theorem}
\numberwithin{equation}{section}

\begin{document}
\title[Extension of Operators]{On Extension of Operators in Classes of Finitely Equivalent Banach Spaces.
Part 3}
\author{Eugene Tokarev}

\section{Bounded extension of operators. 1. A separable case.}

The isomorphic analogue of the class $\frak{P}_{1}(X^{f})$ is the class
$\frak{P}^{\sim}(X^{f})$, which consists of all such spaces $E<_{f}X$, whose
\textit{isomorphic} image in every space from $X^{f}$ admits a bounded projection.

This class essentially differs from $\frak{P}_{1}(X^{f})$ but has much common
with $\frak{J}(X^{f})$. Below it will be shown that (in difference to the
classical case $\frak{P}_{1}(c_{0}^{f})$) there are such classes $X^{f}$ that
$\frak{P}^{\sim}(X^{f})=\varnothing$ and such classes $Y^{f}$ that
$\frak{P}^{\sim}(Y^{f})\neq\varnothing$ but $\frak{P}_{1}(X^{f})\nsubseteq
\frak{P}^{\sim}(X^{f})$.

It is more natural to study the (formally more wide) class $\frak{P}^{\sim
}(X^{F})$, define below.

\begin{definition}
Let $X$ be a Banach space, $X^{F}$ be the corresponding class of crudely
finite equivalence. A space $X_{0}$ is said to be boundedly universally
complemented in $X^{F}$ provided every space $Z\in X^{F}$, which contains a
subspace isometric to $X_{0}$, admits a bounded projection $P:Z\rightarrow
X_{0}$.
\end{definition}

Surely, this definition is equivalent to the formally stronger property of
existence of a bounded projection $P^{\prime}:Z\rightarrow iX_{0}$ where
$i:X_{0}\rightarrow Z$ is an isomorphic embedding.\

Certainly, every finite-dimensional Banach space $A$ generates the class
$A^{F}$, which consists of all finite-dimensional Banach spaces $B$ of
dimension $\dim B=\dim A=n$.

Since for every finite-dimensional Banach space $A$ and every space $Y$ that
contains a subspace $A$ there exists a projection $P:Y\rightarrow A$ of norm
$\left\|  P\right\|  \leq\sqrt{\dim A}$ (cf. [25]), all finite-dimensional
Banach spaces are universally complemented in each class $X^{F}$, which is
generated by an infinite-dimensional space $X.$ So, such a case will be
eliminated as trivial.

A class of all infinite-dimensional Banach spaces, which are boundedly
universally complemented in the class $X^{F}$ will be denoted by
$\frak{P}^{\sim}(X^{F})$.

Similarly to the classes $\mathcal{A}(X^{f})$, $\mathcal{A}_{0}(X^{f})$ and
$\mathcal{E}(X^{f})$ there may be defined their isomorphic analogues --
classes $\mathcal{A}^{\sim}(X^{F})$, $\mathcal{A}_{0}^{\sim}(X^{F})$ and
$\mathcal{E}^{\sim}(X^{F})$.

\begin{definition}
The class $\mathcal{A}^{\sim}(X^{F})$ consists of those
(infinitely-dimensional) Banach spaces $X_{0}$ that have the property:

\begin{itemize}
\item  If $Z\in X^{F}$and $X_{0}\hookrightarrow Z$ then the annihilator
$\left(  X_{0}\right)  ^{\perp}\hookrightarrow Z^{\ast}$ is a complemented
subspace of $Z^{\ast}$.
\end{itemize}

The class $\mathcal{A}_{0}^{\sim}(X^{F})$ is given by $\mathcal{A}_{0}^{\sim
}(X^{F})=\mathcal{A}^{\sim}(X^{F})\cap X^{F}$.
\end{definition}

\begin{definition}
Let $Z\in\mathcal{B}$; $W\hookrightarrow Z$; $\lambda<\infty$. $W$ is said to
be a $\lambda$-reflecting subspace of $Z$ if for every finitely dimensional
subspace $A\hookrightarrow Z$ there exists an isomorphic embedding
$u:A\rightarrow W$ with $\left\|  u\right\|  \left\|  u^{-1}\right\|
\leq\lambda$, such that $u\mid_{A\cap Z}=Id_{A\cap Z}$.
\end{definition}

Certainly if $W$ is a $\lambda$-reflecting subspace of $Z$ for all $\lambda>1$
then $W\prec_{u}Z$.

\begin{definition}
The class $\mathcal{E}^{\sim}(X^{F})$ consists of those
(infinitely-dimensional) Banach spaces $W$ that have the property:

\begin{itemize}
\item  If $Z\in X^{F}$and $W$ is isomorphic to a subspace $iW$ of $Z$ \ then
there exists a such $\lambda<\infty$ that the image $iW$ is a $\lambda
$-reflecting subspace of $Z$.
\end{itemize}
\end{definition}

The following theorem is proved similarly to the isometric analogue. Here
$\frak{qd}^{\sim}$ denotes the class of all Banach spaces $X$ that have a
complement in the second dual $X^{\ast\ast}$ (under the canonical embedding
$k_{X}:X\rightarrow X^{\ast\ast}$).

\begin{theorem}
$\frak{P}^{\sim}(X^{F})=\mathcal{A}^{\sim}(X^{F})\cap\frak{qd}^{\sim}$.

If $Y\in\mathcal{A}^{\sim}(X^{F})$ then $Y^{\ast\ast}\in\frak{P}^{\sim}%
(X^{F})$.

If $X$ is superreflexive then $\frak{P}^{\sim}(X^{F})=\mathcal{A}^{\sim}%
(X^{F})$.
\end{theorem}

This theorem is of restricted interest.

Indeed, in difference with the isometric case, there no analogue to the
non-emptity of the class $\mathcal{E}^{\sim}(X^{F})\subseteq\mathcal{A}^{\sim
}(X^{F})$. Moreover, in some cases even the widest class $\mathcal{A}%
_{0}^{\sim}(X^{F})$ also may be empty.

\begin{example}
If $1\leq p<2$ then $\mathcal{A}_{0}^{\sim}(\left(  l_{p}\right)
^{F})=\mathcal{A}^{\sim}(\left(  l_{p}\right)  ^{F})=\varnothing$.

If $2<p<\infty$ then $\mathcal{A}_{0}^{\sim}(\left(  l_{p}\right)
^{F})=\varnothing$; $\mathcal{A}^{\sim}(\left(  l_{p}\right)  ^{F})=\left(
l_{2}\right)  ^{F}$.
\end{example}

\begin{proof}
Notice that every space $L_{p}$ ($1\leq p\neq2<\infty$) contains an
uncomplemented space, isomorphic to $l_{p}$. Moreover, if $1<p<2$ then $L_{p}$
contains an uncomplemented subspace isomorphic to $l_{2}$ (certainly, this is
obvious for $p=1$). For $p>2$ this result follows from [26] as it was noted in
[27], p. 145. For $1\leq p<4/3$ this is a consequence of one W. Rudin's result
[28]; on the whole interval $\left(  1,2\right)  $ it was extended in [29]. At
last, the existence in $L_{1}$ of an uncomplemented subspace isomorphic to
$l_{1}$ was shown by J. Bourgain [30].

If $W\in\mathcal{A}^{\sim}(\left(  l_{p}\right)  ^{F})$ then it is either
$\mathcal{L}_{p}$-space or $\mathcal{L}_{2}$-space (in the sense of [31]) and,
hence, contains a complemented subspace that is isomorphic to $l_{p}$ (in the
case $\mathcal{L}_{p}$; $1\leq p<\infty$). It is easy to see that the only
possibility is: $p>2$ and $W\in\left(  l_{2}\right)  ^{F}$. Indeed, in this
case any space, which is isomorphic to the Hilbert space, has a complement in
every space of the class $\left(  l_{p}\right)  ^{F}$ that it contains, as it
follows from [32].

In any other case there exists an isomorphic embedding $u$ of $W$ in the space
$L_{p}\left(  \mu\right)  $ (for the suitable measure $\mu$) such that the
image $uW$ has no complement in $L_{p}\left(  \mu\right)  $.
\end{proof}

The result of M.I. Kadec and A. Pelczynski [32], mentioned above, was extended
by B. Maurey [33] to the wide class of Banach spaces of type 2.,

In [33] it was shown:

\begin{itemize}
\item \textit{If a Banach space }$X$\textit{ is of type 2 and }$W$\textit{ is
its subspace, isomorphic to the Hilbert space }$H$,\textit{ then there exists
a projection }$P:X\rightarrow W$\textit{ of norm }%
\[
\left\|  P\right\|  \leq f\left(  T_{p}\left(  X\right)  \right)  d\left(
W,H\right)  ,
\]
\textit{where the constant }$f\left(  T_{p}\left(  X\right)  \right)
$\textit{ depends only on }$T_{p}\left(  X\right)  $\textit{. }
\end{itemize}

An obvious consequence of this result is the following.

\begin{theorem}
For any Banach space $X$ of type 2 the class $\frak{P}^{\sim}(X^{F})$ is non-empty.
\end{theorem}

\begin{proof}
As it follows from the aforementioned B. Maurey's theorem [33], the following
inclusion is fulfilled: $\left(  l_{2}\right)  ^{F}\subseteq\frak{P}^{\sim
}(X^{F})$.
\end{proof}

Notice that in the general (isometric) case classes $\mathcal{A}_{0}(X^{f})$
and $\mathcal{E}(X^{f})$ are different. E.g., if $X=l_{p}$ and $p\neq2$ then
the unique separable space in the class $\mathcal{E}(\left(  l_{p}\right)
^{f})$ is the space $L_{p}\left[  0,1\right]  $. However, spaces $l_{p}$,
$l_{p}\oplus_{p}L_{p}\left[  0,1\right]  $ and a sequence of spaces
$\{l_{p}^{\left(  n\right)  }\oplus_{p}L_{p}\left[  0,1\right]  :n<\infty\}$
are pairwice non isometric and each of them belongs to $\mathcal{A}%
_{0}(\left(  l_{p}\right)  ^{f})$. Notice, also, that in this case
$\mathcal{A}_{0}(\left(  l_{p}\right)  ^{f})=\mathcal{A}(\left(  l_{p}\right)
^{f})$ for all $p\in\left[  1,\infty\right]  $ because of $L_{p}$ contains an
orthogonal complemented subspace $l_{2}$ only if $p=2$..

\begin{theorem}
For every Banach space $X$ classes $\mathcal{E}^{\sim}(X^{F})\ $and
$\mathcal{A}^{\sim}(X^{F})$ are identical.
\end{theorem}

\begin{proof}
The inclusion $\mathcal{E}^{\sim}(X^{F})\subseteq\mathcal{A}^{\sim}(X^{F})$
was noted before. To show the converse recall the following result, due to J.
Stern [9]:

\begin{itemize}
\item \textit{Let }$W$\textit{, }$Z$\textit{ be a pair of Banach spaces;
}$W\hookrightarrow Z$\textit{. If }$W$\textit{ is an algebraic subspace of
}$Z$\textit{ then }$W$\textit{ is a }$\lambda$\textit{-reflecting subspace of
}$Z$\textit{ for some }$\lambda\leq24$\textit{.}
\end{itemize}

Let $Z\in X^{F}$; $Y\in\mathcal{A}^{\sim}(X^{F})$. Without loss of generality
it may be assumed that $Y\in$ $\frak{qd}^{\sim}$ and $Y\hookrightarrow Z$.
Hence, there exists such $\lambda<\infty$ that there exists a projection
$P:Z\rightarrow Y$, $\left\|  P\right\|  \leq\lambda$.

Define on $Z$ a new norm $\left\|  \cdot\right\|  ^{\prime}$ by taking as its
unit ball the convex hull of the union of $\lambda^{-1}B\left(  Y\right)  $
and $B\left(  Z\right)  $ (recall that $B\left(  X\right)  $ denotes the unit
ball of a Banach space $X$):%
\[
B\left(  \left\langle Z,\left\|  \cdot\right\|  ^{\prime}\right\rangle
\right)  =\operatorname{conv}\{(\lambda^{-1}B\left(  Y\right)  )\cup B\left(
Z\right)  \}.
\]
Clearly, $\left\|  y\right\|  ^{\prime}=\left\|  y\right\|  $ for $y\in Y$;
$\left\|  z\right\|  ^{\prime}\leq\lambda\left\|  z\right\|  \leq$
$\lambda\left\|  z\right\|  ^{\prime}\ \ $for $z\in Z$ and there exists a
projection of norm one from $Z^{\prime}=\left\langle Z,\left\|  \cdot\right\|
^{\prime}\right\rangle $ onto $Y.$

Hence $Z^{\prime}$ may be represented as a direct orthogonal sum $Z^{\prime
}=Y\oplus W$, where $W$ is $C$-finitely representable in $Y$ for some
$C<\infty$.

Let $W$ be represented as a direct limit of the direct systems $\left(
W_{\alpha}\right)  _{\alpha\in A}$ of all its finite-dimensional subspaces:
$W=\underset{\rightarrow}{\lim}W_{\alpha}$. Proceeding by induction on the
directed set $A$, chose a set $\left(  W_{\alpha}^{\prime}\right)  _{\alpha\in
A}$ so that it also forms an isometric direct system, spaces $Y\oplus
W_{\alpha}^{\prime}$ are finitely representable in $Y$ and $d(W_{\alpha
},W_{\alpha}^{\prime})\leq C^{\prime}$ for some $C^{\prime}<\infty$. Then the
space $Z^{\prime\prime}=$ $\underset{\rightarrow}{\lim}Y\oplus W_{\alpha
}^{\prime}$ is finitely equivalent to $Y$, isomorphic to $Z$ and $Y$ is an
algebraic subspace of $Z^{\prime\prime}$. By the aforementioned Stern's
theorem, $Y$\ is a $\lambda$-reflecting subspace of $Z^{\prime\prime}$. Using
the inverse isomorphisms we see that $Y$ is $\lambda^{\prime}$-reflecting
subspace of $Z$. Since $Y\in\mathcal{A}^{\sim}(X^{F})$ and $Z\in X^{F}$ are
arbitrary, this yields that $\mathcal{A}^{\sim}(X^{F})\subseteq\mathcal{E}%
^{\sim}(X^{F})$.
\end{proof}

Now it will be shown that the class $\frak{P}^{\sim}(X^{F})$ may be non-empty
for a wide class of Banach spaces that are not of type 2. It will be needed
the following result.

\begin{theorem}
Let a superreflexive class $X^{f}$ contains a separable space $E_{X}$ of
almost universal disposition (with respect to $\frak{M}(X^{f})$). Let $Z_{1}$
and $Z_{2}$ be isomorphic subspaces of $E_{X}$ of equal codimension (in
$E_{X}$); $u:Z_{1}\rightarrow Z_{2}$ be the corresponding isomorphism. Then
there exists an isomorphic automorphism $\widetilde{u}:\ E_{X}\rightarrow
E_{X}$, which restriction to $Z_{1}$ is equal to $u$.
\end{theorem}

\begin{proof}
Let $Z_{1}$ and $Z_{2}$ be subspaces of $E_{X}$ of equal codimension in
$E_{X}$, i.e.
\[
\operatorname{codim}\left(  Z_{1}\right)  =\dim(E_{X}/Z_{1})=\dim(E_{X}%
/Z_{2})=\operatorname{codim}\left(  Z_{2}\right)  .
\]

If $\operatorname{codim}\left(  Z_{1}\right)  <\infty$ then assertions of the
theorem are satisfied independently of a property of $E_{X}$ to be a space of
almost universal disposition. Indeed, for an arbitrary Banach space $X$ any
its subspaces of equal codimension, say, $Y$ and $Z$, are isomorphic and each
of them has a complement in $X$. So, $X$ may be represented both as a direct
sum $X=Y\oplus A$ and as $X=Z\oplus B$, where $\dim(A)=\dim(B)<\infty$.

Let $v:A\rightarrow B$ and $u:Y\rightarrow Z$ be isomorphisms. Clearly,
$v\oplus u:X\rightarrow X$ is an automorphism of $X$, which extends $u$.

Assume that $Z_{1}$ and $Z_{2}$ are isomorphic subspaces of $E_{X}$ of
infinite codimension. Let $u:Z_{1}\rightarrow Z_{2}$ be the corresponding isomorphism.

Let $Z_{1}$ be represented as a direct limit (= as the closure of union) of a
chain
\[
Z_{1}^{\prime}\hookrightarrow Z_{2}^{\prime}\hookrightarrow...\hookrightarrow
Z_{n}^{\prime}\hookrightarrow...\hookrightarrow Z_{1};
\]
namely, $Z_{1}=\overline{\cup\{Z_{n}^{\prime}:n<\omega\}}$ and $\cup
Z_{n}^{\prime}$ be dense in $Z_{1}$. Clearly,
\[
uZ_{1}^{\prime}\hookrightarrow uZ_{2}^{\prime}\hookrightarrow
..\hookrightarrow uZ_{n}^{\prime}\hookrightarrow...\hookrightarrow
uZ_{1}=Z_{2}%
\]
and the union $\cup uZ_{n}^{\prime}$ is dense in $Z_{2}$.

Let $\left(  e_{n}\right)  _{n<\omega}$ be a sequence of linearly independent
elements of $E_{X}$ of norm one, which linear span is dense in $E_{X}$.

Denote the restriction $u\mid_{Z_{n}^{\prime}}$ by $u_{n}$ and define by
induction two sequences $\left(  f_{n}\right)  $ and $\left(  g_{n}\right)  $
of elements of $E_{X}$ and a sequence $\left(  v_{n}\right)  $ of isomorphisms.

Let $f_{1}$ be an element of $\left(  e_{n}\right)  _{n<\omega}$ with the
minimal number, which does not belong to $Z_{1}^{\prime}$.

Put $W_{1}=\operatorname{span}\left(  Z_{1}^{\prime}\cup\{f_{1}\}\right)  $
(recall that $\operatorname{span}(A)$ denotes the \textit{closure} of the
linear span of $A$).

Let $\varepsilon>0$ and $v_{1}:W_{1}\rightarrow E_{X}$ be an extension of
$u_{1}$ such that
\[
\left\|  v_{1}\right\|  \left\|  v_{1}^{-1}\right\|  \leq\left(
1+\varepsilon\right)  \left\|  u_{1}\right\|  \left\|  u_{1}^{-1}\right\|  .
\]
Put $g_{1}=v_{1}\circ f_{1}$; $U_{1}=v_{1}W_{1}$.

Let $g_{2}$ be an element of $\left(  e_{n}\right)  _{n<\omega}$ with the
minimal number, which does not belong to $\operatorname{span}\left(
uZ_{1}^{\prime}\cup U_{1}\right)  $.

Let $U_{2}=\operatorname{span}\left(  uZ_{1}^{\prime}\cup U_{1}\cup
\{g_{2}\}\right)  $.

Since $E_{X}$ is a space of almost universal disposition, there exists an
isomorphism $\left(  v_{2}\right)  ^{-1}$ which extends both $\left(
u_{2}\right)  ^{-1}$ and $\left(  v_{1}\right)  ^{-1}$ and satisfies the
inequality
\[
\left\|  v_{2}\right\|  \left\|  v_{2}^{-1}\right\|  \leq\left(
1+\varepsilon^{2}\right)  \max\{\left\|  u_{2}\right\|  \left\|  u_{2}%
^{-1}\right\|  ;\left\|  v_{1}\right\|  \left\|  v_{1}^{-1}\right\|  \}.
\]
Let $f_{2}=\left(  v_{2}\right)  ^{-1}\circ g_{2}$. This close the first step
of the induction.

Assume that $\{f_{1},f_{2},...,f_{n}\}$; $\{g_{1},g_{2},...,g_{n}\}$;
$\{v_{1},v_{2},...,v_{n}\}$; $\{U_{1},U_{2},...,U_{n}\}$ and $\{W_{1}%
,W_{2},...,W_{n}\}$ are already chosen. If $n$ is odd, we choose sequentially:

$f_{n+1}$to be an element of $\left(  e_{n}\right)  _{n<\omega}$ with the
minimal number, which does not belongs to $Z_{n+1}^{\prime}$;

$W_{n+1}=\operatorname{span}\left(  Z_{n+1}^{\prime}\cup\{f_{n+1}\}\right)  $;

$v_{n+1}:W_{n+1}\rightarrow E_{X}$ be an extension of $u_{n+1}$ such that
\[
\left\|  v_{n+1}\right\|  \left\|  v_{n+1}^{-1}\right\|  \leq\left(
1+\varepsilon^{n+1}\right)  \left\|  u_{n+1}\right\|  \left\|  u_{n+1}%
^{-1}\right\|  ;
\]

$g_{n+1}=v_{n+1}\circ f_{n+1}$; $U_{n+1}=v_{n+1}W_{n+1}$.

If $n$ is even, we choose sequentially:

$g_{n+1}$ to be an element of $\left(  e_{n}\right)  _{n<\omega}$ with the
minimal number, which does not belong to $\operatorname{span}\left(
uZ_{n}^{\prime}\cup U_{n}\right)  $;

$U_{n+1}=\operatorname{span}\left(  uZ_{n}^{\prime}\cup U_{n}\cup
\{g_{n+1}\}\right)  $;

$\left(  v_{n+1}\right)  ^{-1}$, which extends both $\left(  u_{n+1}\right)
^{-1}$ and $\left(  v_{n}\right)  ^{-1}$ and satisfies
\[
\left\|  v_{n+1}\right\|  \left\|  v_{n+1}^{-1}\right\|  \leq\left(
1+\varepsilon^{2}\right)  \max\{\left\|  u_{n+1}\right\|  \left\|
u_{n+1}^{-1}\right\|  ;\left\|  v_{n}\right\|  \left\|  v_{n}^{-1}\right\|
\};
\]

$f_{n+1}=\left(  v_{n+1}\right)  ^{-1}\circ g_{n+1}$.

Since $E_{X}$ is superreflexive, a sequence of isomorphisms $\left(
v_{n}\right)  $ converges to an automorphism $V:E_{X}\rightarrow E_{X}$, which
satisfies
\[
\left\|  V\right\|  \left\|  V^{-1}\right\|  \leq\prod\nolimits_{n=1}^{\infty
}\left(  1+\varepsilon^{n}\right)  \left\|  u\right\|  \left\|  u^{-1}%
\right\|  .
\]
\end{proof}

\begin{theorem}
Let $X^{f}$ be a superreflexive class that contains a separable space $E_{X}$
of almost universal disposition Then $E_{X}\in\mathcal{E}^{\sim}%
(X^{F})=\frak{P}^{\sim}(X^{F})$.
\end{theorem}

\begin{proof}
Let $Z\in X^{F}$; $j:E_{X}\rightarrow Z$ be an isomorphic embedding. Without
loss of generality it may be assumed that $Z$ is separable (if $Z$ is non
separable, it has a separable complemented subspace $Z_{0}$, which contains
$jE_{X}$). Also, it may be assumed that there exists an isometric embedding
$k:Z\rightarrow E_{X}$ (since $E_{X}\in\mathcal{E}(X^{f})$). Hence, the
composition $k\circ j$ is an isomorphic embedding of $E_{X}$ into $E_{X}$.

It may be assumed without the loss of generality that there exists an
isometric embedding $w$ of $E_{X}$ into $E_{X}$.

There exists a bounded projection $P:E_{X}\rightarrow wE_{X}$. According to
the previous theorem, there exists an automorphism $v:E_{X}\rightarrow E_{X}$,
which sends $jE_{X}$ to $wE_{X}$. Clearly, this implies that there exists a
bounded projection from $Z$ onto $jE_{X}$. Since $Z$ is arbitrary, $E_{X}%
\in\frak{P}^{\sim}(X^{F})$.
\end{proof}

Analogously to $\frak{P}^{\sim}(X^{F})$ it may be defined the class
$\frak{J}^{\sim}(X^{F})$.

\begin{definition}
Let $X\ $be a Banach space and$\ E<_{F}X$. $E$ belongs to the class
$\frak{J}^{\sim}(X^{F})$ provided for any pair $Y_{0}\hookrightarrow Y$ where
$Y<_{F}X$ every operator $u:Y_{0}\rightarrow Y$ has a bounded extension
$\widetilde{u}:Y\rightarrow E$.
\end{definition}

\begin{theorem}
Let $X^{f}$ be superreflexive and contains a separable space $E_{X}$ of almost
universal disposition. Let $W$ be isomorphic to $E_{X}$. Then $W\in
\frak{J}^{\sim}(X^{F})$.
\end{theorem}

\begin{proof}
Let $W$ be isomorphic to $E_{X}$; $u:E_{X}\rightarrow W$ be the corresponding
isomorphism. Let $v:W\rightarrow E_{X}$ be the isomorphic embedding.

Let $Y_{0}\hookrightarrow Y$; $Y<_{F}E_{X}$ and $Y$ be separable. Let
$i:Y_{0}\rightarrow Y$ be an operator. Then
\[
v\circ i:Y_{o}\rightarrow vW\hookrightarrow E_{X}%
\]
and this mapping may be extended to the operator $V:Y\rightarrow E_{X}$ of the
same norm.

Put $A=VY\cap vW$.

Since $E_{X}$ is a space of almost universal disposition, there exists a
subspace $W_{0}\hookrightarrow E_{X}$ of infinite codimension such that
$A\hookrightarrow W_{0}$ and $W_{0}\in\mathcal{E}(\left(  E_{X}\right)  ^{F}%
)$. Hence, there exists an automorphism $U:E_{X}\rightarrow E_{X}$ such that
$UW_{0}=W$. Certainly, this isomorphism may be chosen to fix $\left(  v\circ
i\right)  Y_{o}\hookrightarrow E_{X}$.

Clearly, $U\circ V$ extends $v\circ i$ and, hence, $v^{-1}\circ U\circ V$
extends $i$.
\end{proof}

\begin{corollary}
Let $X^{f}$ be a superreflexive class that contains a separable space $E_{X}$
of almost universal disposition. Let $W$ be isomorphic to $E_{X}$. Then
$W\in\frak{P}^{\sim}(X^{F})$.
\end{corollary}

\begin{proof}
As it was shown before, condition on $X^{f}$ implies that $\frak{P}%
(X^{f})=\frak{J}(X^{f})$. Hence, by the previous theorem, $W\in\frak{J}^{\sim
}(X^{F})$. Let $Z\in X_{F}$ and $W\hookrightarrow Z$. Then this identical
embedding may be extended to the operator $P:Z\rightarrow Z$ such that
$P\mid_{W}=Id_{W}$. Clearly, $P$ is a projection.
\end{proof}

\section{Bounded extension of operators. 2. A non-separable case.}

Because of every class of crudely finite equivalence $X^{F}$, which contains a
quotient-closed divisible class $Y^{f}$, contains also its Gurarii compression
$\Gamma(X^{f})$, and, as it was shown before, every superreflexive class of
kind $\Gamma(X^{f})$ contains a separable space of almost universal
disposition $E_{X}$ with respect to $\frak{M}(\Gamma(X^{f}))$, it follows that
for every quotient-closed superreflexive divisible class $X^{f}$ both classes
$\frak{P}^{\sim}(X^{F})$ and $\frak{J}^{\sim}(X^{F})$ are equal and non-empty.
However, it was proved the existence of only one member of these classes,
namely - $E_{X}$. The existence of other spaces without assuming set
theoretical hypothesis is doubtful.

Here it will be considered a non-separable case and it will be shown that
under certain set-theoretical assumptions (e.g., under assumption of existence
an inaccessible cardinal) classes $\frak{P}^{\sim}(X^{F})$ and $\frak{J}%
^{\sim}(X^{F})$ are confinal with $X^{F}$ and are non-empty in the
superreflexive case.

Recall some definitions on cardinal numbers.

\textit{Ordinals} will be denoted by small Greece letters $\alpha,\beta
,\gamma$. \textit{Cardinals} are identified with the least ordinals of given
cardinality and are denoted by $\iota,\tau,\varkappa\ $(may be, with indices).
As usual, $\omega$ and $\omega_{1}$ denote respectively the first infinite and
the first uncountable cardinals (= ordinals).

A set of elements (of arbitrary nature) $\{a_{\alpha}:\alpha<\varkappa\}$ will
be called \textit{the }$\varkappa$\textit{-sequence}.

For a cardinal $\tau$ its \textit{predecessor} (i.e., the least cardinal
$\varkappa$, which is strongly greater then $\tau$) is denoted by $\tau^{+}$.

\textit{The} \textit{confinality of} $\tau,$ $\operatorname{cf}(\tau)$ is the
least cardinality of a set $A\subset\tau$ such that $\tau=\sup A$.

Let $A$, $B$ be sets. The symbol $^{B}A$ denotes the set of all functions from
$B$ to $A$.

In a general case, the cardinality of the set $^{B}A$ is denoted either by
$\operatorname{card}(A)^{\operatorname{card}(B)}$ or by $\varkappa^{\tau}$, if
$\operatorname{card}(A)=\varkappa$; $\operatorname{card}(B)=\tau$.

The symbol $\exp(\tau)$ (or, equivalently, $2^{\tau}$) denotes the cardinality
of the set $\operatorname*{Pow}\left(  \tau\right)  $ of all subsets of $\tau$.

A cardinal $\tau$ is said to be \textit{regular} (resp.,\textit{ singular}) if
its confinality $\operatorname*{cf}\left(  \tau\right)  =\tau$ (resp., if
$\operatorname*{cf}\left(  \tau\right)  <\tau$ ).

A cardinal $\iota$ is said to be \textit{inaccessible} provided it is regular
and $\exp(\varkappa)\leq\iota$ for every $\varkappa<\iota$.

The \textit{continuum hypothesis} (CH) is the assumption%
\[
\omega^{+}\overset{\operatorname*{def}}{=}\omega_{1}=\exp\left(
\omega\right)  .
\]
\

The \textit{general continuum hypothesis} (GCH) is the assumption%
\[
\tau^{+}=\exp\left(  \tau\right)  \text{ \ \textit{for every infinite
cardinal} }\tau.
\]

Let $\tau$ be a cardinal. Put%
\[
\varkappa^{\ast}=\sum\{\exp\left(  \tau\right)  :\tau<\varkappa\}
\]
It is known that there are arbitrary large cardinals $\varkappa$ having the
property
\[
\varkappa=\varkappa^{\ast}.
\]
In other words, a class of all cardinals with the property $\varkappa
=\varkappa^{\ast}$ is confinal with the class of all cardinals. However, there
known such models of the set theory, in which every cardinal $\varkappa$ with
the property $\varkappa=\varkappa^{\ast}$ is singular.

From the other hand, in assuming the general continuum hypothesis (GCH) every
regular cardinal $\varkappa$ has the property $\varkappa=\varkappa^{\ast}$.
Certainly, any inaccessible cardinal $\iota$ is regular and enjoys the
property $\iota=\iota^{\ast}.$

Recall that the existence of an inaccessible cardinal (as well as the CH or
the GCH) cannot be proved in ZFC.

\textbf{In what follows it will be assumed that there exists a regular
cardinal }$\varkappa$\textbf{ with the property }$\varkappa=\varkappa^{\ast}$.
Every such cardinal will be called \textbf{the star-cardinal}.

Recall that a Banach space $X$ is said to be an \textit{envelope} (of the
class $X^{f}$) if every Banach space $Y$, which is finitely representable in
$X$ and is of dimension $\dim Y\leq\dim X$ is isometric to a subspace of $X$.
In [34] it was shown that for every cardinal $\varkappa$ with the property
$\varkappa=\varkappa^{\ast}$ (it does not matter regular or singular) every
class $X^{f}$ of finite equivalence (generated by an infinite-dimensional
space $X$) contains an envelope of this class of dimension $\varkappa$.

It will be convenient to introduce one more definition.

\begin{definition}
A Banach space $X$ is said to be $f$-saturated if $X$ is a space of almost
universal disposition with respect to a class $(X^{<f})_{\dim X}$ of all
Banach spaces that are finitely representable in $X$ and whose dimension is
strictly less then $\dim X$;%
\[
(X^{<f})_{\dim X}=\{Y<_{f}X:\dim Y<\dim X\}.
\]
\end{definition}

\begin{theorem}
Let $X^{f}$ enjoys the isomorphic amalgamation property. Then for every
\textbf{star-cardinal} $\varkappa$ the class $X^{f}$ contains an $f$-saturated
space $E_{\varkappa}\in X^{f}$ of dimension $\varkappa$.
\end{theorem}

\begin{proof}
Let $F$ be an envelope of $X^{f}$ of dimension $\varkappa$. Let $(f_{i}%
)_{i<\varkappa}$ be dense in the unit sphere of $F$.

For every $A\subset\varkappa$ of cardinality $\operatorname*{card}A<\varkappa$
put $F_{A}=\operatorname*{span}\{f_{i}:i\in A\}$. Surely, there are just
$\varkappa$ different spaces of kind $F_{A}$.

For $A$, $B\subset\varkappa$, $\max\{\operatorname*{card}%
A,\operatorname*{card}B\}<\varkappa$ let $u_{A,B}:F_{A}\rightarrow F_{B}$ be
an isomorphic embedding (if it exists). The set of all such embeddings is of
cardinality $\tau\leq\varkappa$.

Let $\varepsilon<1$; $n\in\mathbb{N}$; $i_{A,n}:F_{A}\rightarrow F$ be an
isomorphic embedding of norm $\left\|  i_{A,n}\right\|  \left\|  i_{A,n}%
^{-1}\right\|  \leq1+\varepsilon^{n}$. The set of all such embeddings is of
cardinality that does not exceed $\varkappa$ as well.

So, the set $T$ of all $V$-formations of kind $t=\left\langle F_{A}%
,F,F_{B},i_{A,n},u_{A,B}\right\rangle $ may be indexed by elements of
$\varkappa$. Let $T=\{t_{\alpha}:\alpha<\varkappa\}$.

Construct by induction an $\varkappa$-sequence of spaces $(E_{\alpha}^{\prime
})_{\alpha<\varkappa}$ as a chain
\[
E_{0}^{\prime}\hookrightarrow E_{1}^{\prime}\hookrightarrow...E_{\alpha
}^{\prime}\hookrightarrow....
\]

Put $E_{0}^{\prime}=F$. If $\alpha=\beta+1$ then as $E_{\alpha}^{\prime}$ it
will be chosen an amalgam of the $V$-formation $\left\langle F_{A},E_{\beta
}^{\prime},F_{B},v_{A,n},u_{A,B}\right\rangle $, where $\left\langle
F_{A},F,F_{B},i_{A,n},u_{A,B}\right\rangle $ is just $t_{\beta}$ in our
indexation and $v_{A,n}=j\circ i_{A,n}$, where $j:F\rightarrow E_{\beta
}^{\prime}$ is a natural embedding.

If $\alpha$ is a limit ordinal (i.e. $\alpha$ has no representation
$\alpha=\beta+1$ for any $\beta$) put $E_{\alpha}^{\prime}%
=\operatorname*{span}\{E_{\beta}^{\prime}:\beta<\alpha\}$.

Put $E_{\left(  1\right)  }=\operatorname*{span}\{E_{\alpha}^{\prime}%
:\alpha<\varkappa\}$.

Now we repeat this inductive procedure $\varkappa$ times starting at the
$\gamma$-th step with $E_{\left(  \gamma\right)  }$. Sequentially there will
be obtained $E_{\left(  1\right)  }\hookrightarrow E_{\left(  2\right)
}\hookrightarrow...$.

Put $E_{\varkappa}=\operatorname*{span}\{E_{\left(  \gamma\right)  }%
:\gamma<\varkappa\}$.

Certainly, $\dim E_{\varkappa}=\varkappa$. To show that $E_{\varkappa}$ is
$f$-saturated chose a pair $G$, $H$ of subspaces of $E_{\varkappa}$, $\dim
G\leq\dim H<\varkappa$ and an isomorphic embedding $u:G\rightarrow H$. Since
$F$ is an envelope, $G$ and $H$ almost isometric to corresponding subspaces of
$F$ (say, $F_{G}$ and $F_{H}$). Let $i:G\rightarrow E_{\varkappa}$ be an
isomorphic embedding, $\varepsilon>0$ be fixed.

Since $\varkappa$ is regular, there exists such $\gamma<\varkappa$ that $i$is
$\left(  1+\varepsilon\right)  $-embedding of $G$ into $E_{\left(
\gamma\right)  }$ and, hence, it may be extended to the embedding
$\widetilde{u}$ of $B$ into $E_{\left(  \gamma\right)  }$ with $\left\|
\widetilde{u}\right\|  \left\|  \widetilde{u}^{-1}\right\|  \leq\left(
1+\varepsilon\right)  \left\|  u\right\|  \left\|  u^{-1}\right\|  $.

Since $G$, $H$ and $u:G\rightarrow H$ are arbitrary, this shown that
$E_{\varkappa}$ is $f$-saturated.
\end{proof}

According to [35] a Banach space $X$ is said to be
\textit{subspace-homogeneous} if for every pair $A$, $B$ of its isomorphic
subspaces of equal (finite or infinite codimension) the isomorphism
$v:A\rightarrow B$ may be extended to an isomorphical automorphism
$\widetilde{v}$ of $X$ (i.e., $\widetilde{v}\mid_{A}=v$). Using this
terminology it may be said that the theorem 33 proves that every
superreflexive separable space $E$ of almost universal disposition (with
respect to $\frak{M}(E^{f})$) is subspace-homogeneous. This theorem, that was
proved firstly in the author's preprint [36], solves the problem of [35] on
existence subspace-homogeneous Banach spaces that are differ from $l_{2}$ and
$c_{0}$.

The preceding theorem makes possible to prove the existence of non-separable
subspace-homogeneous spaces (in assumption on existence of star-cardinals).
The question on existence of non-separable subspace-homogeneous spaces also
was posed in [35].

\begin{corollary}
Let $E_{\varkappa}$ be a dual $f$-saturated space of the class $X^{f}$ (this
means that $E_{\varkappa}=W^{\ast}$ for some $W\in\mathcal{B}$). Then
$E_{\varkappa}$ is subspace-homogeneous.

In particular, $E_{\varkappa}\in\frak{P}^{\sim}(X^{F})$ (= $\frak{J}^{\sim
}(X^{F})$)
\end{corollary}

\begin{proof}
The proof literally repeats the proof of theorems 33 and 34 with obvious
changes: in the present proof the transfinite induction is used.. At the last
step of the proof is used the $w^{\ast}$-compactness of the unit ball of
$E_{\varkappa}$.
\end{proof}

\begin{theorem}
Let $E$ be $f$-saturated and dual (i.e., $E=W^{\ast}$ for some $W\in
\mathcal{B}$; in particular, $E$ may be separable); $F$ be a complemented
subspace of $E$. Then $F\in\frak{P}^{\sim}(X^{F})$ (= $\frak{J}^{\sim}(X^{F})$).
\end{theorem}

\begin{proof}
Let $A$, $B$ be isomorphic subspaces of $E$; $A\hookrightarrow
F\hookrightarrow E$; $B\hookrightarrow F\hookrightarrow E$, which are of equal
(finite or infinite) codimension in $F$. Let $u:A\rightarrow B$ be the
corresponding isomorphism. Let $U:E\rightarrow E$ be an authomorphism of $E$
that extends $u$; $Ua=ua$ for all $a\in A$.

Let $P:E\rightarrow F$ be a projection. Certainly, $UF\hookrightarrow E$ is
also a complemented subspace of $E$. Let $Q:E\rightarrow UF$ be the
corresponding projection

So, $F\cap UF$ is also a complemented subspace of $E$, which contains $B$.
Clearly, $F\cap UF$ is a complemented subspace of both $F$ and $UF$. Let
$F^{\prime}$ be its complement in $F$ (resp., let $F^{\prime\prime}$ be its
complement in $UF$). Clearly, $F^{\prime}$ and $F^{\prime\prime}$ are
isomorphic. Let $v:F^{\prime}\rightarrow F^{\prime\prime}$ be the
corresponding isomorphism; $id$ be identical on $F\cap UF$. Then the operator
$v\circ id\circ U\mid_{F}is$an automorphism of $F$, which restriction to $A$
is equal to $u$.
\end{proof}

\begin{remark}
In $\left[  36\right]  $ it was shown that if $E$ is subspace-homogeneous and
its $l_{p}$-spectrum $S(E)$ contains a point $p\neq2$ then $E$ and $E\oplus E$
are non-isomorphic. So, from the preceding theorem it follows that there
exists a pair of non-isomorphic separable spaces of $\frak{P}^{\sim}(X^{F})$
(= $\frak{J}^{\sim}(X^{F})$) when $E$ satisfies the mentioned conditions.
\end{remark}

\begin{corollary}
In assumption the GCH every superreflexive quotient-closed divisible class
$X^{f}$ contains spaces from $\frak{P}^{\sim}(X^{F})$ (and, thus, from
$\frak{J}^{\sim}(X^{F})$) of arbitrary large cardinality.

The same is true if we assume the existence of inaccesible cardinals.
\end{corollary}

Now it will be shown that the condition of regularity of $\varkappa$ cannot be
omitted, at least, in the case when the class $X^{f}$ is not super-stable in
the sense of [37-38]. Recall the definition.

\begin{definition}
($\left[  37\right]  $). A Banach space $X$ is said to be stable provided for
any two sequences $\left(  x_{n}\right)  $ and $\left(  y_{m}\right)  $ of its
elements and every pair of ultrafilter $D$, $E$ over $\mathbb{N}$
\[
\lim_{D\left(  n\right)  }\lim_{E\left(  m\right)  }\left\|  x_{n}%
+y_{m}\right\|  =\lim_{E\left(  m\right)  }\lim_{D\left(  n\right)  }\left\|
x_{n}+y_{m}\right\|  .
\]
\end{definition}

The notations $D(n)$ and $E(m)$ are used here (instead of $D$ and $E$) to
underline the variable ($n$ or $m$ respectively) in expressions like
$\lim_{D\left(  n\right)  }f(n,m)$.

\begin{definition}
($[38]$). A Banach space $X$ is said to be superstable if every its ultrapower
$\left(  X\right)  _{D}$ is stable.
\end{definition}

Let $X$ be a Banach space.

A sequence $\{x_{n}:n<\infty\}$ of elements of $X$ is said to be

\begin{itemize}
\item \textit{Spreading,} if for any $n<\infty$, any $\varepsilon>0$, any
scalars $\{a_{k}:k<n\}$ and any choosing of $i_{0}<i_{1}<...<i_{n-1}<...$;
$j_{0}<j_{1}<...<j_{n-1}<...$ of natural numbers
\[
\left\|  \sum\nolimits_{k<n}a_{k}x_{i_{k}}\right\|  =\left\|  \sum
\nolimits_{k<n}a_{k}x_{j_{k}}\right\|  .
\]

\item \textit{Symmetric,} if for any $n<\omega$, any finite subset
$I\subset\mathbb{N}$ of cardinality $n$, any rearrangement $\varsigma$ of
elements of $I$ and any scalars $\{a_{i}:i\in I\}$,
\[
\left\|  \sum\nolimits_{i\in I}a_{i}z_{i}\right\|  =\left\|  \sum
\nolimits_{i\in I}a_{\varsigma(i)}z_{i}\right\|  .
\]
\end{itemize}

\begin{definition}
Let $X$ be a Banach space. Its $IS$-spectrum $IS(X)$ is a set of all
(separable) spaces $\left\langle Y,\left(  y_{i}\right)  \right\rangle $ with
a spreading basis $\left(  y_{i}\right)  $ which are finitely representable in
$X$.
\end{definition}

\begin{theorem}
A class $X^{f}$ is superstable if and only if every member $\left\langle
Y,\left(  y_{i}\right)  \right\rangle $ of its $IS$-spectrum has a symmetric basis.
\end{theorem}

\begin{proof}
The first part of the theorem was proved in [37].

Conversely, let every $\left\langle Z,\left(  z_{i}\right)  \right\rangle $
has a symmetric basis. Suppose that $X$ is not superstable. Then there exists
a space from $X^{f}$ which is not stable (it may be assumed that $X$ is not
stable itself). By [37] there are such sequences $\left(  x_{n}\right)  $ and
$\left(  y_{m}\right)  $ of elements of $X$ that
\[
\sup_{m<n}\left\|  x_{n}+y_{m}\right\|  >\inf_{m>n}\left\|  x_{n}%
+y_{m}\right\|  .
\]

Let $D$ be a countably incomplete ultrafilter over $\mathbb{N}$. Put
\[
X_{0}\overset{\operatorname{def}}{=}X\text{; \ }X_{n}\overset
{\operatorname{def}}{=}\left(  X_{n-1}\right)  _{D};\text{ }n=1,\text{
}2,\text{ }...;\text{ \ }X_{\infty}=\overline{\cup_{n\geq1}X_{n}}.
\]
Here is assumed that $X_{n}$ is a subspace of $X_{n+1}=\left(  X_{n}\right)
_{D}$ under the canonical embedding $d_{X_{n}}:X_{n}\rightarrow\left(
X_{n}\right)  _{D}$.

Let $D$, $E$ be ultrafilters over $\mathbb{N}$. Their product $D\times E$ is a
set of all subsets $A$ of $\mathbb{N}\times\mathbb{N}$ that are given by
\[
\{j\in\mathbb{N}:\{i\in\mathbb{N}:\left(  i,j\right)  \in A\}\in D\}\in E.
\]
Certainly, $D\times E$ is an ultrafilter and for every Banach space $Z$ the
ultrapower $\left(  Z\right)  _{D\times E}$ may be in a natural way identified
with $\left(  \left(  Z\right)  _{D}\right)  _{E}$.

So, the sequence $\left(  x_{n}\right)  \subset X$ defines elements
\begin{align*}
\frak{x}_{1}  &  =\left(  x_{n}\right)  _{D}\in\left(  X\right)  _{D};\\
\frak{x}_{2}  &  =\left(  x_{n}\right)  _{D\times D}\in\left(  \left(
X\right)  _{D}\right)  _{D};\\
&  ...\text{ \ \ \ \ }...\text{ \ \ \ \ }...\text{ \ \ \ \ }...\text{
\ \ \ \ }...\text{ \ \ \ \ }...\\
\frak{x}_{k}  &  =\left(  x_{n}\right)  \underset{k\text{ }times}%
{_{\underbrace{D\times D\times...\times D}}}\in\underset{k\text{ }%
times}{\underbrace{(\left(  \left(  X\right)  _{D}\right)  _{D}...)_{D}};}\\
&  ...\text{ \ \ \ \ }...\text{ \ \ \ \ }...\text{ \ \ \ \ }...\text{
\ \ \ \ }...\text{ \ \ \ \ }...
\end{align*}

Notice that $\frak{x}_{k}\in X_{k}\backslash X_{k-1}$. It is easy to verify
that $\left(  \frak{x}_{k}\right)  _{k<\infty}\subset X_{\infty}$ is a
spreading sequence. Since $X_{\infty}\in X^{f}$, it is symmetric. Moreover,
for any $z\in X$, where $X$ is regarded as a subspace of $X_{\infty}$ under
the direct limit of compositions
\[
d_{X_{n}}\circ d_{X_{n-1}}\circ...\circ d_{X_{0}}:X\rightarrow X_{n},
\]
the following equality is satisfied: for any pair $m$, $n\in\mathbb{N}$
\[
\left\|  x_{n}+z\right\|  =\left\|  x_{m}+z\right\|  .
\]
Since $\left(  x_{n}\right)  $ and $z$ are arbitrary elements of $X$, this
contradicts with the inequality $\sup_{m<n}\left\|  x_{n}+y_{m}\right\|
>\inf_{m>n}\left\|  x_{n}+y_{m}\right\|  $.
\end{proof}

\begin{theorem}
Let $E$ be an envelope of the class $E^{f}$, which is not super-stable. Let
$\dim E=\varkappa$ and $\varkappa$ is singular. Then $E$ is not $f$-saturated.
\end{theorem}

\begin{proof}
Since $E^{f}$ is not super-stable, there exists a space $\langle
W,(w_{n})\rangle$ with a spreading basis, which is not symmetric one. Consider
a vector space $c_{00}\left(  \varkappa\right)  $ of all $\varkappa$-sequences
of reals all but finitely many members of which are vanished.

Let $e_{\alpha}=(\delta_{\alpha\beta})_{\beta<\varkappa}\in c_{00}\left(
\varkappa\right)  $ ($\alpha<\varkappa$), where $\delta_{\alpha\beta}$ be the
generalized Kronecker symbol; $\delta_{\alpha\beta}=0$ when $\alpha\neq\beta$;
$\delta_{\alpha\beta}=1$ when $\alpha=\beta$.

Every element $x\in c_{00}\left(  \varkappa\right)  $ is of kind $x_{A}%
=\sum\{x_{\alpha}e_{\alpha}:\alpha\in A\}$, where $A\subset\varkappa$ is finite.

Let $A=\{\alpha_{0},\alpha_{1},...,\alpha_{m}\}$ and $\alpha_{0}<\alpha
_{1}<...<\alpha_{m}$. Define on $c_{00}\left(  \varkappa\right)  $ a norm,
which is given by%
\[
\left\|  x_{A}\right\|  =\left\|  \sum\nolimits_{i=0}^{m}x_{\alpha_{i}}%
w_{i}\right\|  _{W}.
\]
Let $W(\varkappa)$ be a completition of $c_{00}\left(  \varkappa\right)  $
under this norm.

Surely, $W(\varkappa)$ has a spreading, non symmetric basis, is finitely
representable in $E_{\varkappa}$ and is isometric to a subspace of
$E_{\varkappa}$ (since $\dim W(\varkappa)=\varkappa=\dim E_{\varkappa}$ and
$E_{\varkappa}$is an envelope).

Since $\varkappa$ is singular, $\operatorname*{cf}(\varkappa)=\tau<\varkappa$.
Let $\left(  \varkappa_{\gamma}\right)  _{\gamma<\tau}$ be an increasing
$\tau$-sequence of cardinals (that are less then $\varkappa$), such that
$\sup\{\varkappa_{\gamma}:\gamma<\tau\}$

Consider two subspaces of $W(\varkappa)$, say $W_{1}(\tau)$ and $W_{2}(\tau)$.

$W_{1}(\tau)$ is spaned by the first $\tau$ elements of $e_{\alpha}$;%
\[
W_{1}(\tau)=\operatorname*{span}\{e_{\alpha}:\alpha<\tau\}.
\]

$W_{2}(\tau)$ is spanned by those $e_{\alpha}$'s, whose indices belong to the
set $\left(  \varkappa_{\gamma}\right)  _{\gamma<\tau}$;%
\[
W_{1}(\tau)=\operatorname*{span}\{e_{\tau_{\beta}}:\beta<\tau\}.
\]

Surely, $W_{1}(\tau)$ and $W_{2}(\tau)$ are isometric, but no automorphism
$u:E\rightarrow E$ sends $W_{1}(\tau)$ to $W_{2}(\tau)$.
\end{proof}

\section{References}

\begin{enumerate}
\item  Nachbin I. \textit{A theorem of Hahn-Banach type for linear
transformations}, Trans. AMS \textbf{68} (1950) 28-46

\item  Goodner D. A. \textit{Projections in normed and linear spaces}, Trans.
AMS \textbf{69} (1950) 89-108

\item  Kelley J. \textit{Banach spaces with the extension property}, Trans.
AMS \textbf{72} (1952) 323-326

\item  Lindenstrauss J. \textit{On a problem of Nachbin concerning extension
of operators}, Israel J. Math. \textbf{2:1} (1963) 75-90

\item  Lindenstrauss J. \textit{Extension of compact operators}, Memoirs AMS
\textbf{48} (1964)

\item  Tokarev E.V. \textit{Injective Banach spaces in the finite equivalence
classes }(transl. from Russian), Ukrainian Mathematical Journal \textbf{39:6}
(1987) 614-619

\item  Dacunha-Castelle D., Krivine J. -L. \textit{Applications des
ultraproduits \`{a} l'\'{e}tude des espaces et des alg\`{e}bres de Banach},
Studia Math. \textbf{41} (1972) 315 - 334

\item  Schwartz L. \textit{Geometry and probability in Banach spaces}, Bull.
AMS \textbf{4:2} (1981) 135-141

\item  Stern J. \textit{Ultrapowers and local properties in Banach spaces},
Trans. AMS \textbf{240} (1978) 231-252

\item  Yasuhara M. \textit{The amalgamation property, the universal
homogeneous models and generic models}, Math. Scand. \textbf{34} (1974) 5 - 36

\item  Lindenstrauss J., Rosenthal H.P.\textit{ The} $\mathcal{L}_{p}%
$\textit{-spaces}, Israel J. Math.\textbf{ 7} (1969) 325-349

\item  Davis W.J., Figiel T., Johnson W.B., Pe\l czy\'{n}ski A.
\textit{Factoring weakly compact operators}, J. Funct. Anal. \textbf{17}
(1974) 311-327

\item  Heinrich S. \textit{Ultraproducts in Banach space theory}, J. Reine
Angew. Math. \textbf{313} (1980) 72-104

\item  Gurarii V.I. \textit{On inclinations and spreadings of subspaces of
Banach spaces} (in Russian), Teor. Funct., Funct. Analysis and Appl.\textbf{1}
(1965) 194-204

\item  Pisier G. \textit{On the duality between type and cotype}, Lect. Notes
in Math. \textbf{939} (1982)

\item  Pietsch A. \textit{Absolut }$p$\textit{-summierende Abbildungen in
normierten R\"{a}umen}, Studia Math. \textbf{28} (1967) 333-353

\item  Pe\l czy\'{n}ski A. $p$\textit{-integral operators commuting with group
representations and examples of quasi }$p$\textit{-integral operators, which
are not }$p$\textit{-integral}, Studia Math. \textbf{33} (1969) 63-70

\item  Pisier G. \textit{Counterexample to a conjecture of Grothendieck}, Acta
Math. \textbf{151} (1983) 181-208

\item  Grothendieck A. \textit{Sur quelques points d'algebre homologique},
T\^{o}hoku Math. J., second series, \textbf{9:2,3} (1957) 119-221

\item  Heinrich S. \textit{Finite representability and super - ideals of
operators}, Diss. Math. (Rozprawy Matematyczne) \textbf{172} (1980) Warszawa

\item  Gurarii V.I.\textit{\ Spaces of universal disposition, isotropic spaces
and the Mazur problem on rotations in Banach spaces}, Sibirsk. Mat. Journ. (in
Russian) \textbf{7} (1966) 1002-1013

\item  Maurey B., Pisier G. \textit{S\'{e}ries de variables al\'{e}atoires
vectorielles ind\'{e}pendantes et propri\'{e}t\'{e}s g\'{e}om\'{e}triques des
espaces de Banach,} Studia Math. \textbf{58} (1976) 45-90

\item  Rosenthal H.P. \textit{On a theorem of J. -L. Krivine concerning block
finite representation of }$l_{p}$\textit{ in general Banach spaces}, J. Funct.
Anal. \textbf{25} (1978) 197-225

\item  Stromquist W. \textit{The maximum distance between two-dimensional
spaces}, Math. Scand. \textbf{48} (1981) 205-225

\item  Kadec M.I., Snobar L.G. \textit{On some functionals on Minkowski's
spaces} (in Russian), Mat. Zametki \textbf{10} (1971) 453-458

\item  Rosenthal H.P. \textit{Projections onto translation-invariant subspaces
of }$L_{p}(G)$, Memoirs AMS \textbf{63} (1968) 1-84

\item  Lindenstrauss J. and Tzafriri L. \textit{Classical Banach spaces,}
Lecture Notes in Math. \textbf{338} (1973) 1-242

\item  Rudin W. \textit{Trigonometric series with gaps}, J. Math. Mech.
\textbf{9} (1960) 203-227

\item  Bennet G., Dor L.E., Gudman V., Johnson W.B., Newmann C.M. \textit{On
uncomplemented subspaces of }$L_{p}$, $1<p<2$, Israel J. Math. \textbf{26}
(1977) 178-187

\item  Bourgain J. \textit{Complementation de sous-espaces }$L_{1}$\textit{
dans les espaces }$L_{1}$, Sem. Funct. Anal. \'{E}cole Polytech. Palaiseau,
exp. \textbf{27} (1980) 1-7

\item  Lindenstrauss J., Pe\l czy\'{n}ski A. \textit{Absolutely summing
operators in }$\mathcal{L}_{p}$\textit{-spaces and their applications}, Studia
Math. \textbf{29 }(1968) 275-326

\item  Kadec M.I. and Pe\l czy\'{n}ski A. \textit{Bases, lacunary sequences
and complemented subspaces in the spaces }$L_{p}$ Studia Math. \textbf{21}
(1962) 161-176

\item  Maurey B. \textit{Quelques probl\`{e}mes de factorization
d'op\'{e}rateurs lin\'{e}aires}, Actes du Congress Int. Math. Vancouver, 1974
v. \textbf{2} (1975) 75-79

\item  Tokarev E.V. \textit{Problem of envelopes of locally convex spaces}
(transl. from Russian), Siberian Mathematical Journal \textbf{31:1} (1990) 173-175

\item  Lindenstrauss J., Rosenthal H.P. \textit{Authomorphisms in }$c_{0}%
$\textit{, }$\ell_{1}$\textit{ and }$m$, Israel J. Math. \textbf{7} (1969) 227-239

\item  Tokarev E.V. \textit{A solution of one problem of Lindenstrauss and
Rosenthal on subspace-homogeneous and quotient-homogeneous Banach spaces with
application to the approximation problem and to the Schroeder - Bernstein
problem}, Preprint. Printed electronically in \textbf{www.ArXiv.org} \ at
June, 03, 2002,\textbf{ math.FA/0206013}

\item  Krivine J. -L., Maurey B. \textit{Espaces de Banach stables}, Israel J.
Math. \textbf{39} (1981) 273-295

\item  Raynaud Y. \textit{Espaces de Banach superstables}, C. R. Acad. Sci.
Paris, S\'{e}r. A. \textbf{292:14} (1981) 671-673
\end{enumerate}
\end{document}